\mag=\magstep 1										

\documentclass{article}  
\usepackage{romp}
\title{The entropy of co-compact open covers}

\begin{document}

\maketitle
\hspace*{2.90cm} {\large Zheng Wei}$^{\: a}$,$\;\:${\large Yangeng Wang}$^{\: b,1}$,$\;$ {\large Guo Wei}$^{\:c,}$\symbolfootnote[1]{Corresponding author. tel.: (910) 521 6582; fax: (910) 522 5755.\\ \hspace*{0.4cm} {\it E-mail addresses:} weizheng@nmsu.edu (Z. Wei), ygwang@nwu.edu.cn (Y. Wang), guo.wei@uncp.edu (G. Wei)\\ \hspace*{0.4cm}$^1$Supported in Part by SF of Shaanxi (98SL06), P. R. China.}

\bigskip
\hspace*{0.6cm} $^a$ {\it \small Department of Mathematical Science, New Mexico State University, Las Cruces,}\\
\hspace*{4.52cm}                    {\it\small  New Mexico,88001, U.S.A.}\\
\hspace*{0.6cm} $^b$ {\it \small Department of Mathematic, Northwest University, Xian, Shaanxi 710069, P.R.China.}\\
\hspace*{1.6cm} $^c$ {\it \small Department of Mathematics $\&$ Computer Science, University of North Carolina}\\
\hspace*{3.52cm}      {\it \small at Pembroke, Pembroke, North Carolina 28372, U.S.A.}\\


\bigskip

\begin{abstract}
Co-compact entropy is introduced as an invariant of topological conjugation
for perfect mappings defined on any Hausdorff space
(compactness and metrizability not necessarily required). This is achieved through the consideration of
co-compact covers of the space.
The advantages of co-compact entropy include: 1) it does not require the space to be compact, and thus generalizes
Adler, Konheim and McAndrew's topological entropy of continuous mappings on compact dynamical systems, and
2) it is an invariant of topological conjugation, compared to Bowen's entropy that is metric-dependent.
Other properties of co-compact entropy are investigated, e.g., the co-compact entropy of a subsystem does not exceed
that of the whole system.  
For the linear system $(R, f)$ defined by $f(x) = 2x$, the co-compact entropy is zero, while Bowen's entropy for
this system is at least $\log 2$. More general, it is found that co-compact entropy is
a lower bound of Bowen's entropies, and the proof of this result generates the Lebesgue Covering Theorem to co-compact
open covers of non-compact metric spaces, too.
\end{abstract}

\bigskip


\bigskip


\section{\hspace*{-0.15in}. Introduction}
\subsection{Measure-theoretic entropy}         
The concept of entropy per unit time was introduced by Shannon \cite{Shannon}, by analogy with the standard Boltzmann entropy measuring a spatial disorder in a thermodynamic system. In 1950s, Kolmogorov \cite{Kolmogorov} and Sinai established a rigorous definition of K-S entropy per unit time for dynamical systems and other random processes \cite{Cornfeld}. Kolmogorov imported Shannon's probabilistic notion of entropy into the theory of dynamical systems, and the idea was vindicated later by Ornstein who showed that metric entropy suffices to completely classify two-sided Bernoulli processes \cite{Ornstein}, a basic problem which for many decades appeared completely intractable. Kolmogorov's metric entropy is an invariant of measure theoretical dynamical systems and is closely related to Shannon's source entropy. The K-S entropy is a powerful concept because it controls the top of the hierarchy of ergodic properties: K-S property $\Rightarrow$ multiple mixing $\Rightarrow$ mixing $\Rightarrow$ weak mixing $\Rightarrow$ ergodicity \cite{Cornfeld}.
The K-S property holds if there exists a subalgebra of measurable sets in phase space which generates the whole algebra by application of the flow \cite{Cornfeld}. The dynamical randomness of a deterministic system finds its origin in the dynamical instability and the sensitivity to initial conditions. In fact, the K-S entropy is related to the Lyapunov exponents according to a generalization of Pesin's theorem \cite{Eckmann, Kantz}. A deterministic system with a finite number of degrees of freedom is chaotic if its K-S entropy per unit time is positive. More properties about K-S entropy can be found in papers \cite{Cornfeld, Eckmann, van Beijeren}. The concept of space-time entropy or entropy per unit time and unit volume was later introduced by Sinai and Chernov \cite{Sinai}. A spatially extended system with a probability measure being invariant under space and time translations can be said to be chaotic if its space-time entropy is positive.

\subsection{\hspace*{-0.15in}. Topological entropy and its relation to measure-theoretic entropy}
In 1965, Adler, Konheim and McAndrew introduced the concept of topological entropy for continuous mappings defined on compact spaces \cite{Adler}, which is an analogous invariant under conjugation of topological dynamical systems and can be obtained by maximizing the metric entropy over a suitable class of measures defined on a dynamical system, implying that topological entropy and measure-theoretic entropy are closely related. Motivated by a conjecture of Adler, Konheim and McAndrew, Goodwyn in 1969 and 1971 compared topological entropy and measure-theoretic entropy and concluded that topological entropy bounds measure-theoretic entropy \cite{Goodwyn,Goodwyn2}. Bowen in 1970 studied topological entropy and Axiom A \cite{Bowen6}, and generalized the concept of topological entropy to continuous mappings defined on metric spaces and proved that his definition coincides with that of Adler, Konheim and McAndrew's within the class of compact metric spaces. In 1971, Bowen also considered the entropy for non-compact sets and
for group endomorphisms and homogeneous spaces respectively \cite{Bowen5, Bowen2}.
However, the entropy according to Bowen's definition is metric-dependent (\cite{Walters}, Walters' book, p.169) and can be positive even for a linear function (Walters' book, p.176). In 1973, along with a study of measure-theoretic entropy, Bowen in \cite{Bowen5} gave another definition of topological entropy resembling Hausdorff dimension, which also equals to the topological entropy defined by Adler, Konheim and McAndrew when the space is compact. Recently, Liu, Wang and Wei, Canovas and Rodriguez,  Malziri, and Molaci proposed other definitions of topological entropy for continuous mappings defined on non-compact (metric) spaces \cite{Liu, Canovas, Malziri}.

\subsection{\hspace*{-0.15in}. The importance of entropy}
The concept of entropy is useful for studying topological and measure-theoretic structures of dynamical systems. For instance, two conjugated systems have a same entropy and thus entropy is a numerical invariant of the class of conjugated dynamical systems. The theory of expansive dynamical systems has been closely related to the theory of topological entropy \cite{Bowen4, Keynes, Thomas}. Entropy and chaos are closely related, e.g., a continuous mapping $f: I \rightarrow I$ is chaotic if and only if it has a positive topological entropy \cite{Block}. A remarkable result is that a deterministic system together with an invariant probability measure defines a random process. As a consequence, a deterministic system can generate dynamical randomness, which is characterized by an entropy per unit time that measures the disorder of the trajectories along the time axis. Entropy has many applications, e.g., transport properties in escape-rate theory \cite{Gaspard1, Gaspard2, Dorfman, Gaspard3, Gaspard5}, where an escape of trajectories is introduced by absorbing conditions at the boundaries of a system. These absorbing boundary conditions select a set of phase-space trajectories, forming a chaotic and fractal repeller, which is related to an equation for K-S entropy. The escape-rate formalism has applications in diffusion \cite{Gaspard4}, to reaction-diffusion \cite{Claus}, and recently to viscosity \cite{Viscardy}. Another application is the classification of quantum dynamical systems, which is given by Ohya \cite{Ohya}. Symbolic dynamical systems $(\sum(p),\sigma)$ have various representative and complicated dynamical properties and characteristics, with an entropy $\log \:p$. When determining whether or not a given topological dynamical system has certain dynamical complexity, it is often compared with a symbolic dynamical system \cite{Robinson, Zhou}. For the topological conjugation with symbolic dynamical systems, we refer to Ornstein \cite{Ornstein}, Sinai \cite{Sinai1}, Akashi \cite{Akashi}, and Wang and Wei  \cite{Wang, Wang2}.

\subsection{\hspace*{-0.15in}. The purpose, the approach and the outlines}
The main purpose is to introduce a topological entropy for perfect mappings defined on arbitrary Hausdorff spaces (compactness and metrizability not necessarily required), and investigate fundamental properties of such an entropy.

Instead of using all open covers of the space to define entropy, we consider the open covers consisting of the co-compact
open sets (open sets whose complements are compact).

Various definitions of entropy and historical notes are mentioned previously in this section.
Section 2 investigates the topological properties of co-compact open covers of a space.
Section 3 introduces the new topological entropy defined through co-compact covers of the space,
which is called co-compact entropy in the paper.
Section 4 further explores the properties of the co-compact entropy and compares it with Adler, Konheim and McAndrew's topological
entropy for compact spaces.
Sections 5 and 6 investigate the relation between the co-compact entropy and Bowen's entropy. More precisely,
Section 5 compares the co-compact  entropy with that given by Bowen for systems defined on metric spaces.
Because the spaces under consideration include non-compact metric spaces, the traditional Lebesgue Covering
Theorem does not apply. Thus, one work is to generalize this theorem
to co-compact open covers of non-compact metric spaces.
Based on the generalized Lebesgue Covering Theorem, we show that the co-compact entropy is a lower bound for
Bowen's entropies. In Section 6, a linear dynamical system is studied. For this simple system, its co-compact entropy is $0$ which is appropriate,
but Bowen's entropy is positive.

\section{\hspace*{-0.15in}. Basic concepts and definitions}			
Let $(X, f)$ be a topological dynamical system, where $X$ is a Hausdorff and $f : X \rightarrow X$ is a continuous mapping.
We introduce the concept of co-compact open covers as follows.

\begin{definition} {\rm Let $X$ be a Hausdorff space. For an open subset $U$ of $X$, if $X \backslash U$ is a compact subset of $X$,
then $U$ is called a co-compact open subset. If every element of an open cover $\mathcal{U}$ of $X$  is co-compact,
then $\mathcal {U}$ is called a co-compact open cover of $X$. }
\end{definition}

\begin{theorem}  \label{theorem2.1}
The meet of finitely many co-compact open subsets is co-compact, and the union of any collection of co-compact open subsets is
co-compact open.
\end{theorem}

\noindent
{\bf Proof.}
Suppose that $U_{1},U_{2},...,U_{n}$ are co-compact open. Let $U = \bigcap\limits_{i=1}^{n}U_{i}$.
As $X \setminus U_{i}, i=1,2,...n$, are compact, $X \setminus U = \bigcup\limits_{i=1}^{n} (X  \setminus U_{i})$
is compact and hence $U$ is co-compact open.

Suppose that $\{U_{\lambda}\}_{\lambda\in\Lambda}$ is a family of co-compact sets. Let
$U=\bigcup\limits_{\lambda\in\Lambda} U_\lambda$. As for any $\lambda\in\Lambda$ $X \backslash U_\lambda$ is
compact, $X \backslash U = \bigcap\limits_{\lambda\in\Lambda} (X\backslash U_\lambda)$ is compact. Hence,
$U$ is co-compact open.
$\diamondsuit$

\begin{theorem}\label{theorem2.2} Let $X$ be Hausdorff. Then any co-compact open cover has a finite subcover.
\end{theorem}

{\bf Proof.} Let $\mathcal{U}$ be a co-compact open cover. For any $U \in \mathcal {U}$, $X \backslash U$ is compact.
Noting that $\mathcal {U}$ is also an open cover of $X \setminus U$, there exists a finite subcover $\mathcal {V}$ of $X \setminus U$.
Now, $\mathcal {V} \cup \{U\}$ is finite subcover of  $\mathcal {U}$.
$\diamondsuit$

\begin{definition}  {\rm
Let $X$ and $Y$ be  Hausdorff  spaces and $f : X \rightarrow Y$ a continuous mapping. If $f$ is a closed mapping
and all fibers $f^{-1}(x),  x \in Y$, are compact,  then $f$ is called a perfect mapping. }
\end{definition}

In particular, if $X$ is compact Hausdorff and $Y$ is Hausdorff, every continuous mapping from $X$ into $Y$ is perfect.
If $f:X \rightarrow Y$ is perfect, then $f^{-1}(F)$ is compact for each compact subset $F \subseteq Y$ \cite{Engelking}.

\begin{theorem}  \label{theorem2.3}
Let $X$ and $Y$ be two Hausdorff  spaces and $f: X \rightarrow Y$ a perfect mapping. If $U$ is co-compact  open in $Y$, then $f^{-1}(U)$ is
co-compact open in $X$. Moreover,  if $\mathcal{U}$ is a co-compact open cover of $Y$, then  $f^{-1}(\mathcal{U})$ is a co-compact
open cover of $X$.
\end{theorem}

\noindent
{\bf Proof.} It suffices to show that the pre-image of any co-compact set is co-compact. Let $U$ be co-compact open in $Y$.
Then $F = Y \setminus U$ is compact in $Y$. As $f$ is perfect, $f^{-1}(F)$ is compact in $X$.
Hence, $f^{-1}(U) = X \setminus f^{-1}(F)$ is co-compact open in $X$.
$\diamondsuit$

\section{\hspace*{-0.15in}. The Entropy of co-compact open covers}
For compact topological systems, Adler, Konheim and McAndrew introduced the concept of topological entropy and studied its
properties  \cite{Adler}. Their definition is as follows: Let $X$ be a compact topological space and $f: X \rightarrow X$ a continuous
mapping. For any open cover $\mathcal{U}$ of $X$, denote by
$N_X(\mathcal{U})$ the smallest cardinality of all subcovers of $\mathcal{U}$, i.e.,
$$N_X(\mathcal{U})=\min\{card(\mathcal{V}) : \mathcal{V} {\rm \; is\; a \; subcover\; of\;} \mathcal{U}\}.$$
It is obvious that  $N_X(\mathcal{U})$ is a positive integer. Let $H_X(\mathcal{U})=\log N_X(\mathcal{U})$.
Then $ent(f, \mathcal{U}, X) =
\lim\limits_{n\rightarrow\infty}\frac{1}{n}H_X(\bigvee\limits_{i=0}^{n-1}f^{-i}(\mathcal{U}))$
is called the topological entropy of $f$  relative to
$\mathcal{U}$, and  $ent(f,X)=\sup\limits_\mathcal{U}\{ent(f,\mathcal{U},X)\}$ is called the topological entropy of  $f$.

Now, we will generalize Adler, Konheim and McAndrew's entropy to any Hausdorff space for perfect mappings.
So in the remainder of the paper, a space is assumed to be Hausdorff and a mapping is assumed to be perfect.

Let $X$ be Hausdorff. By Theorem \ref{theorem2.2}, when $\mathcal {U}$ is a co-compact open cover of $X$, $\mathcal {U}$
has a finite subcover. Hence, $N_X(\mathcal {U})$, abbreviated as $N(\mathcal {U})$, is a positive integer. Let
$H_X(\mathcal {U}) = \log \: N(\mathcal {U})$, abbreviated as $H(\mathcal {U})$.

Let $\mathcal {U}$ and $\mathcal {V}$ be two open covers of $X$. Define

$$\mathcal{U}\bigvee\mathcal{V}=\{U \cap V : U\in\mathcal{U} {\rm \; and\;} V\in\mathcal{V}\}.$$

If for any $U \in \mathcal {U}$, there exists $V \in \mathcal {V}$ such that $U \subseteq V$, then $\mathcal {U}$
is said to be a refinement of $\mathcal {V}$ and is denoted by $\mathcal{V}\prec\mathcal{U}$.

The following are some obvious facts:

\smallskip
\smallskip

{\bf Fact 1:} For any open covers $\mathcal{U}$ and $\mathcal{V}$ of $X$,  $\mathcal{U}\prec\mathcal{U}\bigvee\mathcal{V}$.
\smallskip
\smallskip

{\bf Fact 2:} For any open covers $\mathcal{U}$ and $\mathcal {V}$ of $X$,  if $\mathcal{V}$ is a subcover of $\mathcal{U}$, then $\mathcal{U}\prec\mathcal{V}$.
\smallskip
\smallskip

{\bf Fact 3:} For any co-compact open cover $\mathcal{U}$ of $X$,
$H(\mathcal{U})=0 \Longleftrightarrow N(\mathcal{U})=1 \Longleftrightarrow X\in\mathcal{U}$.
\smallskip
\smallskip

{\bf Fact 4:}  For any co-compact open covers $\mathcal{U}$ and $\mathcal {V}$ of $X$,
$\mathcal{V}\prec\mathcal{U}\Rightarrow H(\mathcal{V}) \leq H(\mathcal{U})$.
\smallskip
\smallskip

{\bf  Fact 5:} For any co-compact open covers $\mathcal{U}$ and $\mathcal{V}$,
$H(\mathcal{U}\bigvee\mathcal{V})\leq H(\mathcal{U})+H(\mathcal{V})$.

To prove Fact 5, let $\mathcal {U}_0$ be a finite subcover of $\mathcal {U}$, with the cardinality $N(\mathcal {U})$.
Let $\mathcal {V}_0$ be a finite subcover of $\mathcal {V}$ with the cardinality $H(\mathcal {V})$. Then
$\mathcal{U}_{0}\bigvee\mathcal{V}_{0}$ is a subcover of $\mathcal{U}\bigvee\mathcal{V}$, and the cardinality of $\mathcal{U}_{0}\bigvee\mathcal{V}_{0}$ is at most
$N(\mathcal{U}) \times N(\mathcal{V})$. Hence,
$N(\mathcal{U}\bigvee\mathcal{V}) \leq N(\mathcal{U}) \times N(\mathcal{V})$, and therefore
$H(\mathcal{U}\bigvee\mathcal{V}) \leq H(\mathcal{U})+H(\mathcal{V})$.

\smallskip
\smallskip

{\bf Fact 6:} For any co-compact open cover $\mathcal{U}$ of $X$, $H(f^{-1}(\mathcal{U}))\leq H(\mathcal{U})$, and if $f(X) = X$
the equality holds.

To prove Fact 6, let  $\mathcal {U}_0$ be a finite subcover of $\mathcal {U}$, with the cardinality $N(\mathcal {U})$.
$f^{-1}(\mathcal {U}_0)$ is a subcover of $f^{-1}(\mathcal {U})$. Hence, we have $H(f^{-1}(\mathcal {U}) \leq H(\mathcal {U})$.

Now, assume $f(X)=X$. Let $\{f^{-1}(U_{1}),f^{-1}(U_{2}),...,f^{-1}(U_{n})\}$,
$U_{i}\in\mathcal{U}$, be a finite subcover of $f^{-1}(\mathcal{U})$, with the cardinality $N(f^{-1}(\mathcal{U}))$. As
$X \subseteq \bigcup\limits_{i=1}^{n}f^{-1}(U_{i})$, we have
$X=f(X)\subseteq \bigcup\limits_{i=1}^{n}f(f^{-1}(U_{i}))=\bigcup\limits_{i=1}^{n}U_{i}$. Hence,
${U_{1},U_{2},...,U_{n}}$ is a finite subcover of $\mathcal {U}$. This shows $H(\mathcal {U})\leq H(f^{-1}(\mathcal {U}))$. This
inequality and the previous inequality together imply the required equality.

\begin{lemma} \label{lemma3.1}
Let $\{a_{n}\}_{n=1}^{ \infty }$ be a sequence of non-negative real numbers satisfying $a_{n+p}\leq a_{n}+a_{p}, n \geq 1, p\geq 1$.
Then $\lim\limits_{n\rightarrow\infty}\frac {a_{n}}{n}$ exists and equals to $\inf \frac {a_{n}}{n}$ (see \cite{Walters}).  $\diamondsuit$
\end{lemma}

Let $\mathcal {U}$ be a co-compact open cover of $X$. By Theorem \ref{theorem2.3}, for any positive integer $n$ and perfect mapping
$f : X \rightarrow X$, $f^{-n}(\mathcal {U})$ is a co-compact open cover of $X$. On the other hand, by Theorem \ref{theorem2.1},
$\bigvee\limits_{i=0}^{n-1}f^{-i}(\mathcal{U})$ is a co-compact open cover of $X$. These two facts together lead to
the following result:

\begin{theorem}\label{theorem3.2}
Suppose that  $X$ is Hausdorff. Let $\mathcal{U}$ be a co-compact open cover of $X$, and $f : X \rightarrow X$ a perfect mapping.
Then $\lim\limits_{n\rightarrow\infty}\frac {1}{n}H(\bigvee\limits_{i=0}^{n-1}f^{-i}(\mathcal{U}))$ exists.
\end{theorem}

\noindent
{\bf Proof.}
Let $a_{n}=H(\bigvee\limits_{i=0}^{n-1}f^{-i}(\mathcal{U}))$. By Lemma \ref{lemma3.1},
it suffices to show $a_{n+k}\leq a_{n}+a_{k}$. Now, Fcat 6 gives $H(f^{-1}(\mathcal{U}))\leq H(\mathcal{U})$, and
more general $H(f^{-j}(\mathcal{U}))\leq H(\mathcal{U}), j=0,1,2,...$.
Hence, by applying Fact 5,  we have
$a_{n+k}= H(\bigvee\limits_{i=0}^{n+k-1}f^{-i}(\mathcal{U}))
=H((\bigvee\limits_{i=0}^{n-1}f^{-i}(\mathcal{U}))\bigvee
(\bigvee\limits_{j=n}^{n+k-1}f^{-j}(\mathcal{U})))=
H(\bigvee\limits_{i=0}^{n-1}f^{-i}(\mathcal{U})\bigvee(\bigvee
\limits_{j=0}^{k-1}f^{-n}(f^{-j}(\mathcal{U})))) \leq
H(\bigvee\limits_{i=0}^{n-1}f^{-i}(\mathcal{U}))+
H(\bigvee\limits_{j=0}^{k-1}f^{-j}(\mathcal{U})) =a_{n}+a_{k}$.
$\diamondsuit$

\bigskip

Next, we introduce the concept of entropy for co-compact open covers.

\begin{definition}  {\rm
Let $X$ be a Hausdorff space, $f: X\rightarrow X$ be a perfect mapping, and $\mathcal{U}$ be a co-compact open cover of $X$.
The non-negative number $c(f,\mathcal{U})=\lim\limits_{n \rightarrow \infty} \frac {1}{n} H(\bigvee\limits_{i=0}^{n-1}f^{-i}(\mathcal{U}))$ is said to be
the co-compact entropy of $f$ relative to $\mathcal {U}$, and the non-negative number
$c(f)=\sup\limits_{\mathcal {U}} \{c(f, \mathcal{U}) \}$ is said to be the co-compact entropy of $f$. }
\end{definition}

In particular, when $X$ is compact Hausdorff, any open set of $X$ is co-compact and any continuous mapping $f : X \rightarrow X$
is perfect. Hence, Adler, Konheim and McAndrew's topological entropy is a special case of our co-compact entropy. It should be
aware that the new entropy is well defined for perfect mappings on non-compact spaces, e.g., on $R^n$, but
Adler, Konheim and McAndrew's topological entropy requires that the space be compact.

Co-compact entropy generalizes Adler, Konheim and McAndrew's topological entropy,
and yet it holds various  similar properties as well, as demostrated by the fact that co-compact entropy
is an invariant of topological conjugation (next theorem)  and more explored in the next section.

Recall that $ent$ denotes Adler, Konheim and McAndrew's topological entropy, and $c$
denotes the co-compact entropy.

\begin{theorem} \label{theorem3.4}
Let $(X,f)$ and $(Y,g)$  be two topological dynamical systems where $X$ and $Y$ are Hausdorff, $f : X \rightarrow X$ and
$g : Y \rightarrow $ are perfect mappings.
If there exists a semi-topological conjugation $h: X \rightarrow Y$ where $h$ is also perfect, then $c(f)\geq c(g)$.
Consequently, when  $h$ is a topological conjugation, we have $c(f)=c(g)$.
\end{theorem}

\noindent
{\bf Proof.}
Let $\mathcal{U}$ be any co-compact open cover of $Y$. As
$h$ is perfect and $\mathcal{U}$ is co-compact open cover of $Y$,
$h^{-1}(\mathcal{U})$ is co-compact open cover of $X$ by applying Theorem \ref{theorem2.3}.
Hence, we have

$$c(g,\mathcal{U})=\lim\limits_{n\rightarrow\infty}\frac {1}{n}H(\bigvee\limits_{i=0}^{n-1}g^{-i}(\mathcal{U}))=\lim\limits_{n\rightarrow\infty}\frac {1}{n}H(h^{-1}(\bigvee\limits_{i=0}^{n-1}g^{-i}(\mathcal{U})))$$ $$\hspace*{2cm}=\lim\limits_{n\rightarrow\infty}\frac {1}{n}H(\bigvee\limits_{i=0}^{n-1}h^{-1}(g^{-i}(\mathcal{U})))
=\lim\limits_{n\rightarrow\infty}\frac {1}{n}H(\bigvee\limits_{i=0}^{n-1}f^{-i}(h^{-1}(\mathcal{U})))$$
$$\hspace*{-3.01cm}=c(f,h^{-1}(\mathcal{U}))\leq c(f,\mathcal{U}).$$
Therefore, $c(f) \geq c(g)$.

When $h$ is a topological conjugation, it is of course perfect, too. Hence, we have both $c(f) \geq c(g)$ and
$c(g)\geq c(f)$ from above proof, implying $c(f)=c(g)$.
$\diamondsuit$

\section{\hspace*{-0.15in}. Properties of co-compact entropy}   \label{sec4}			
In this section, we investigate further properties of the co-compact entropy. These properties are comparable to that of
Adler, Konheim and McAndrew's topological entropy.

\begin{theorem} \label{theorem4.1}
Let $X$ be Hausdorff and $id : X\rightarrow X$ be the identity mapping. Then $c(id)=0$.
\end{theorem}

\noindent
{\bf Proof.}
Let $\mathcal{U}$ be any co-compact open cover of $X$.
Then we have $c(id,\mathcal{U})=\lim\limits_{n\rightarrow\infty}\frac {1}{n}H(\bigvee\limits_{i=0}^{n-1}id^{-i}(\mathcal{U}))$
$=\lim\limits_{n\rightarrow\infty}\frac {1}{n}H(\mathcal{U})=0$. Hence,  $c(id)=0$.
$\diamondsuit$

\bigskip

When $X$ is Hausdorff and $f:X\rightarrow X$ is perfect, $f^{m}: X \rightarrow X$ is also a perfect mapping \cite{Engelking}.

\begin{theorem}  \label{theorem4.2}.
Let $X$ be Hausdorff and $f:X\rightarrow X$ be perfect. Then $c(f^{m})=m\cdot c(f)$.
\end{theorem}

\noindent
{\bf Proof.}
Let $\mathcal{U}$ be any co-compact open cover of $X$. As
$\bigvee\limits_{j=0}^{n-1}(f^{m})^{-j}(\bigvee\limits_{i=0}^{m-1}f^{-i}(\mathcal{U}))=
\bigvee\limits_{j=0}^{mn-1}f^{-j}(\mathcal{U})$, we have
$H(\bigvee\limits_{j=0}^{n-1}(f^{m})^{-j}(\bigvee\limits_{i=0}^{m-1}f^{-i}(\mathcal{U})))=
H(\bigvee\limits_{j=0}^{mn-1}f^{-j}(\mathcal{U}))$.
Put $\mathcal{V}=\bigvee\limits_{i=0}^{m-1}f^{-i}(\mathcal{U})$. Then
$c(f^{m})\geq c(f^{m},\mathcal{V})=\lim\limits_{n\rightarrow\infty}\frac
{1}{n}H(\bigvee\limits_{j=0}^{n-1}(f^{m})^{-j}(\bigvee\limits_{i=0}^{m-1}f^{-i}(\mathcal{U}))
=\lim\limits_{n\rightarrow\infty}m  \frac
{1}{mn}H(\bigvee\limits_{j=0}^{mn-1}f^{-j}(\mathcal{U}))=m\cdot c(f,\mathcal{U})$,  thus
$c(f^{m})\geq m\cdot \sup\limits_\mathcal{U}\{c(f,\mathcal{U})\}=m\cdot c(f)$.

On the other hand, it follows from
$$\mathcal{U}\bigvee(f^{m})^{-1}(\mathcal{U})\bigvee...\bigvee (f^{m})^{-(n-1)}(\mathcal{U})\prec \mathcal{U}\bigvee
f^{-1}(\mathcal{U})\bigvee...\bigvee f^{-mn+1}(\mathcal{U})$$ that
$$c(f^{m},\mathcal{U})=\lim\limits_{n\rightarrow\infty}\frac {1}{n}H(\bigvee\limits_{j=0}^{n-1}(f^{m})^{-j}(\mathcal{U}))\leq
m\lim\limits_{n\rightarrow\infty}\frac {1}{mn}H(\bigvee\limits_{j=0}^{mn-1}f^{-j}(\mathcal{U}))=m \cdot c(f,\mathcal{U}).$$
Hence,
$$c(f^{m})=\sup\limits_\mathcal{U}\{c(f^{m},\mathcal{U})\}\leq m\sup\limits_\mathcal{U}\{c(f,\mathcal{U})\}=m \cdot c(f).$$

Therefore, $c(f^{m})=m \cdot c(f)$.
$\diamondsuit$

\begin{theorem} \label{theorem4.3}
Let $X$ be Hausdorff and $f:X\rightarrow X$ be perfect. If $\Lambda$ is a closed subset of $X$ and invariant under $f$, i.e.,
$f(\Lambda) \subseteq \Lambda$, then  $c(f |_{\Lambda}) \leq c(f)$.
\end{theorem}

\noindent
{\bf Proof.}
Let $\Gamma$ denote the collection of all co-compact open cover of $\Lambda$.
For any $\mathcal{U} \in \Gamma$, put $\mathcal{U}^* = \{U \cup (X \setminus \Lambda) \mid U \in \mathcal{U}\}$.
Then $\mathcal{U}^{*}$  is a co-compact open cover of  $X$, and $H(\bigvee\limits_{i=0}^{n-1}(f\mid_{\Lambda})^{-i}(\mathcal{U}))\leq H(\bigvee\limits_{i=0}^{n-1}f^{-i}(\mathcal{U}^{*}))$. Hence, we have
$c(f\mid_{\Lambda})=\sup\limits_\mathcal{U} \{c(f\mid_{\Lambda},\mathcal{U})\}\leq  \sup\limits_\mathcal{U}\{c(f,\mathcal{U}^{*})\}\leq c(f)$.
$\diamondsuit$

\section{\hspace*{-0.15in}. Relations between co-compact entropy and Bowen's  entropy}   \label{sec5}	 
\subsection{\hspace*{-0.15in}. Co-compact entropy less than or equal to Bowen's entropy, $c(f) \leq h_d(f)$}
First let us recall the definition of Bowen's entropy \cite{Bowen6, Walters}. Let $(X, d)$ be a metric space and
$f : X \rightarrow X$ a continuous mapping. A compact subset $E$ of $X$ is called a $(n, \epsilon)$-separated set
with respect to $f$ if for any different $x, y \in E$, there exists an integer $j$ with $0 \leq j < n$ such that
$d(f^j(x), f^j(y)) > \epsilon$. A subset $F$ of $X$ is called a $(n, \epsilon)$-spanning set of a compact set
$K$ relative to $f$ if for any $x \in K$, there exists $y \in F$ such that for all $j$ satisfying $0 \leq j < n$,
$d(f^j(x), f^j(y)) \leq \epsilon$.

Let $K$ be a compact subset of $X$. Put
$$r_{n}(\epsilon,K,f)=\min\{card(F) :  F {\rm \; is\; a\;} (n,\epsilon){\rm -spanning \; set \;for \;} K {\rm \; with\; respect \; to\;} f \},$$
$$s_{n}(\epsilon,K,f)=\max\{card(F) :  F \subseteq K {\rm \; and\;} F {\rm \; is\; a\;} (n,\epsilon) {\rm -separated \; set \; with \; respect \; to \;} f\},$$
$$r(\epsilon,K,f)=\lim\limits_{n\rightarrow\infty}\frac
{1}{n}\log  r_{n}(\epsilon,K,f),~~~~~~~s(\epsilon,K,f)=\lim\limits_{n\rightarrow\infty}\frac {1}{n}\log s_{n}(\epsilon,K,f),$$
$$r(K,f)=\lim\limits_{\epsilon\rightarrow 0} r(\epsilon,K,f) ,~~~~~~~s(K,f)=\lim\limits_{\epsilon\rightarrow 0} s(\epsilon,K,f).$$
Then $\sup\limits_{K}r(K,f)=\sup\limits_{K}s(K,f)$, and this non-negative number denoted by $h_d(f)$ is the Bowen entropy of $f$.

It should be pointed out that Bowen's entropy $h_d(f)$ is metric-dependent, see e.g. \cite{Walters, Liu}. For
the topology of the metrizable space $X$, the selection of different metrics may result in different entropies.

Next, recall the Lebesgue Covering Theorem and Lebesgue Number \cite{Engelking}. Let $(X, d)$ be a metric space and
$\mathcal {U}$ an open cover of $X$. $diam(\mathcal {U}) = \sup \{d(A) \mid A\in \mathcal{U}\}$ is called the
diameter of $\mathcal {U}$, where $d(A)=\sup\{d(x,y)\mid x,y\in A\}$. A real number $\delta$ is said to be
a Lebesgue Number of $\mathcal {U}$ if every open subset $U$ of $X$ satisfying $diam(U) < \delta$ is completely
contained in an element of the cover $\mathcal {U}$.

\bigskip

\noindent
{\bf The Lebesgue Covering Theorem} (see \cite{Engelking}): Every open cover of a compact metric space has a Lebesgue number.
$\diamondsuit$

\bigskip

Our next theorem generalizes the Lebesgue Covering Theorem to all co-compact open covers of
non-compact metric spaces.

\begin{theorem} \label{lemma5.3}
Let $(X,d)$ be a metric space, regardless of compactness. Then every co-compact open cover of $X$ has a Lebesgue number.
\end{theorem}

\noindent
{\bf Proof.}
Let $\mathcal {U}$ be any co-compact open cover of $X$.  By Theorem \ref{theorem2.2}, $\mathcal {U}$ has a finite
subcover $\mathcal{V}=\{V_{1},V_{2},...,V_{m}\}$. Put
$Y=(X \setminus V_{1})\cup (X \setminus V_{2})\cup ...\cup (X \setminus V_{m})$. Then $Y$ is compact as
$V_i$'s are co-compact.

We will prove that $\mathcal {V}$ has a Lebesgue number, so does $\mathcal {U}$. As it is obvious that
the theorem holds when $Y = \emptyset$, thus in the following proof we assume $Y \neq \emptyset$.

Assuming in contradition that $\mathcal {V}$ does not have a Lebesgue number.
Then for any positive integer $n$, $\frac{1}{n}$ is not a Lebesgue number of $\mathcal {V}$. Consequently,
for each positive integer $n$, there exists an open subset $O_n$ of $X$ satisfying $diam(O_n) < \frac{1}{n}$
but $O_n$ is not completely contained
in any element of $\mathcal {V}$, i.e., $O_n \cap (X \setminus V_j) \neq \emptyset, j = 1, 2, ..., m$. Hence,
$O_n \cap Y \neq \emptyset$. Take an $x_n \in O_n \cap Y$.
By the compactness of $Y$, the sequence $x_n$  has a
subsequence $x_{n_i}$ that is convergent to some point $y \in Y$, i.e.,
$\lim\limits_{i \rightarrow\infty}x_{n_{i}}= y \in Y \subseteq  X$.

On the other hand, $\mathcal {V}$ is an open cover of $X$, thus there exists some $V \in \mathcal {V}$ such that
$y \in V$. As $V$ is open, there exists an open neighborhood $S(y, \epsilon)$ of $y$ such that
$y \in S(y, \epsilon) \subseteq V$. Since $x_{n_i}$ converges to $y$, there exists a positive integer
$M$ such that $x_{n_i} \in S(y, \frac{\epsilon}{2})$ for $i > M$. Let $k$ be any integer larger than $M + \frac{2}{\epsilon}$.
Then for any $z \in O_{n_k}$, we have $d(z,y)\leq d(z,x_{n_{k}})+d(x_{n_{k}} , y)< \frac {\epsilon} {2} + \frac
{\epsilon}{2}  =  \epsilon $, thus $O_{n_{k}} \subseteq S(y,\epsilon) \subseteq V \in \mathcal{V}$, which contradicts
the selection of open sets $O_n$'s.

Therefore, $\mathcal {V}$ has a Lebesgue number.
$\diamondsuit$

\begin{theorem} \label{theorem6.4}
Let $(X, d)$ be a metric space, $\mathcal{U}$ be any co-compact open cover of $X$, and $f: X \rightarrow X$ be
a perfect mapping. Then there exists $\delta > 0$ and a compact subset $K$ of $X$ such that for all positive
integers $n$,
$$N(\bigvee\limits_{i=0}^{n-1}f^{-i}(\mathcal{U}))\leq n \cdot r_{n}(\frac{\delta}{3},K,f)+1.$$
\end{theorem}

\noindent
{\bf Proof.}
Let $\mathcal{U}$ be any co-compact open cover of $X$. By Theorem \ref{theorem2.2},  $\mathcal {U}$ has a finite
subcover $\mathcal{V}=\{V_{1},V_{2},...,V_{m}\}$. By Theorem \ref{lemma5.3},
$\mathcal{U}$ has a Lebesgue number $\delta$.
Put $K=(X \setminus V_{1}) \cup(X \setminus V_{2})\cup...\cup(X \setminus V_{m})$.
If $K = \emptyset$, then $X = V_j$ for all $j = 1, 2, ..., m$ and in this case the theorem clearly holds. Hence, we assume
$K \neq \emptyset$, thus  the compact set
$K$ has a $(n,\frac{\delta}{3})$-spanning set $F$ relative to $f$ and satisfying $card(F)=  r_{n}(\frac{\delta}{3},K,f)$.

a) For any $x\in K$, there exists $y\in F$ such that $d(f^{i}(x),f^{i}(y))\leq \frac {\delta}{3}, i=0, 1,...,n-1$, equivalently,
$x \in f^{-i}(\overline{S(f^{i}(y),\frac{\delta}{3})}), i=0,1,...,n-1$.
Hence, $K\subseteq \bigcup\limits_{y\in F}\bigcap\limits_{i=0}^{n-1} f^{-i}(\overline{S(f^{i}(y),\frac{\delta}{3})})$.
By the definition of Lebesgue number, every $\overline{S(f^{i}(y),\frac{\delta}{3})}$ is a subset of an element of $\mathcal {V}$.
Hence, $\bigcap\limits_{i=0}^{n-1}f^{-i}(\overline{S(f^{i}(y),\frac{\delta}{3})})$ is a subset of an element of
$\bigvee\limits_{i=0}^{n-1}f^{-i}(\mathcal{V})$. Consequently, $K$ can be covered by $r_{n}(\frac {\delta}{3},K,f)$
elements of $\bigvee\limits_{i=0}^{n-1}f^{-i}(\mathcal{V})$.

b) For any  $x \in  X \setminus K$, i.e., $x \in V_{1}\cap V_{2}\cap...\cap V_{m}$.
In the following, we will consider points of $X \setminus K$ according to two further types of points.

First, consider those $x$ for which there exists $l$ with $1\leq l\leq n-1$, such that $f^{l}(x) \in K$ and $x,f(x),f^{2}(x), ...,f^{l-1}(x)\in X \setminus K$ ($l$ depends on $x$ but for convenience, we use $l$ instead of $l_x$). Namely,
we consider the set $\{x \in X \setminus K \::\: x \in X \setminus K, x, f(x), f^{2}(x), ...,f^{l-1}(x) \in X  \setminus K, f^l(x) \in K\}$.
For every such $x$, there exists $y\in F$, such that $d(f^{l+i}(x),f^{i}(y))\leq \frac {\delta}{3}, i=0,1,...,n-l-1$, equivalently,
$x\in f^{-(l+i)}(\overline{S(f^{i}(y),\frac{\delta}{3})}),i=0,1,...,n-l-1$.
By the definition of Lebesgue number, every $\overline{S(f^{i}(y),\frac{\delta}{3})}$
is a subset of an element of $\mathcal{V}$. Hence, $V_{1}\cap f^{-1}(V_{1})\cap ... \cap f^{-(l-1)}(V_{1})\cap( \bigcap\limits_{i=0}^{n-l-1}
f^{-(l+i)}(\overline{S(f^{i}(y), \frac{\delta}{3})}))$ is a subset of an element of $\bigvee\limits_{i=0}^{n-1}f^{-i}(\mathcal{V})$, and
$x\in V_{1}\cap f^{-1}(V_{1})\cap ...\cap f^{-(l-1)}(V_{1})\cap( \bigcap\limits_{i=0}^{n-l-1}f^{-(l+i)}(\overline{S(f^{i}(y),\frac{\delta}{3})}))$.
There are $r_{n}(\frac {\delta}{3},K,f)$ such open sets, implying that
$\bigvee\limits_{i=0}^{n-1}f^{-i}(\mathcal{V})$ has $r_{n}(\frac {\delta}{3},K,f)$ elements
that cover this type of points $x$. As $1\leq l\leq n-1$,  $\bigvee\limits_{i=0}^{n-1}f^{-i}(\mathcal{V})$ has $(n-1) \cdot r_{n}(\frac{\delta}{3},K,f)$ elements that actually cover this type of points  $x$.

Next, consider those $x$ for which $f^{i}(x)\in X \setminus K$ for every $i = 0, 1, ..., n-1$.  One (any)  element of
$\bigvee\limits_{i=0}^{n-1}f^{-i}(\mathcal{V})$ covers all such points $x$. Hence,
$X \setminus K$ can be covered by no more than  $(n-1) \cdot r_{n}(\frac {\delta}{3},K,f)+1$ elements
of $\bigvee\limits_{i=0}^{n-1}f^{-i}(\mathcal{V})$.

By a) and b), for any $n > 0$, it holds
$N(\bigvee\limits_{i=0}^{n-1}f^{-i}(\mathcal{V}))\leq n\cdot r_{n}(\frac {\delta}{3},K,f)+1$.
Now, it follows from $\mathcal{U} \prec \mathcal{V}$ and Fact 4, we have
$N(\bigvee\limits_{i=0}^{n-1}f^{-i}(\mathcal{U}))\leq N(\bigvee\limits_{i=0}^{n-1}f^{-i}(\mathcal{V}))\leq n\cdot r_{n}(\frac
{\delta}{3},K,f)+1$.
$\diamondsuit$

\begin{theorem} \label{theorem6.5}
Let $(X, d)$ be a metric space
and $f:X\rightarrow X$ be a perfect mapping. Then $c(f) \leq  h_{d}(f)$.
\end{theorem}

\noindent
{\bf Proof.} For any co-compact open cover $\mathcal {U}$ of $X$, if $X \in \mathcal {U}$, then $c(f, \mathcal {U}) = 0$. Hence,
we can assume $X \not\in \mathcal {U}$.
By Theorem \ref{theorem6.4}, there exists $\delta > 0$ and a non-empty compact subset $K$ of $X$ such that
for any $n > 0$, it holds $N(\bigvee\limits_{i=0}^{n-1}f^{-i}(\mathcal{U}))\leq n \cdot r_{n}(\frac {\delta}{3},K,f)+1$.
Hence, $c(f,\mathcal{U})=\lim\limits_{n\rightarrow\infty}\frac {1}{n}H(\bigvee\limits_{i=0}^{n-1}f^{-i}(\mathcal{U}))$
$\leq \lim\limits_{n\rightarrow\infty}\frac {1}{n}\log(n\cdot r_{n}(\frac
{\delta}{3},K,f)+1)=r(\frac{\delta}{3},K, f)$. Let $\delta\rightarrow 0$. It follows from the definition of Bowen's entropy
(Walters'  book \cite{Walters}, P.168, Definition 7.8 and Remark 2) that $r(\frac{\delta}{3},K, f)$ is decreasing
on $\delta$ and $r(K,f)=\lim\limits_{\delta\rightarrow 0}  r(\frac{\delta}{3},K,f)$. Therefore, we have
$c(f,\mathcal{U})\leq r(\frac{\delta}{3},K, f) \leq r(K,f)$. Moreover, $r(K,f)\leq h_{d}(f)$. Finally, because $\mathcal {U}$
is arbitrarily selected, we have $c(f)\leq h_{d}(f)$.
$\diamondsuit$

\bigskip

Bowen's entropy $h_d(f)$ is metric-dependent. Theorem \ref{theorem6.5} indicates that the co-compact entropy, which is
metric-independent, is always bounded by Bowen's entropy, i.e., $c(f)\leq h_d(f)$, regardless of the choice of a metric for the calculation
of Bowen's entropy. In the next section, we will give an example where co-compact entropy is strictly less than Bowen's entropy.

\subsection{\hspace*{-0.15in}. An example}   \label{SEC5}
In this section, $R$ denotes the one-dimensional Euclidean space equipped with the usual metric $d(x, y) =|x - y|, x, y \in R$.
The mapping $f: R \rightarrow R$ is defined by $f(x)=2x, x \in R$. $f$ is clearly a perfect mapping.
It is known that $h_d(f) \geq \log 2$ \cite{Walters}.  We will show $c(f) = 0$.

Let $\mathcal{V}$ be any co-compact open cover of $R$. By Theorem \ref{theorem2.2},
$\mathcal {V}$ has a finite co-compact subcover $\mathcal {U}$. Let $m = card(\mathcal {U})$.
As compact subsets of $R$ are  closed and bounded sets, there exist
$U_r, U_{l} \in \mathcal{U}$ such that for any $U \in \mathcal{U}$,
$\sup{\{R \setminus U\}} \leq \sup{\{R \setminus U_{r}\}}$ and
$\inf{\{R \setminus U\}}\geq \inf{\{R \setminus U_{l}\}}$. Let $a_{r}=\sup{\{R \setminus U_{r}\}}$
and $b_{l}=\inf{\{R \setminus U_{l}\}}$. Observe that for any
$n > 0$, $x\in\bigvee\limits_{i=0}^{n-1}f^{-i}(U_{i})\Longleftrightarrow x\in U_{0},f(x)\in
U_{1}, ..., f^{n-1}(x)\in U_{n-1}$, where $U_{i}\in \mathcal{U}, i=0,1,...,n-1$.

\smallskip
\smallskip

Case 1: $0<b_{l}<a_{r}$. For any $n > 0$ and $x\in(a_{r},+\infty)$, we have $x\in U_{r},f(x)\in
U_{r},...,f^{n-1}(x)\in U_{r}$. So $(a_{r},+\infty) \subseteq \bigcap\limits_{i=0}^{n-1}f^{-i}(U_{r})$.
For any $x\in(-\infty,0]$, we have $x\in U_{l},f(x)\in U_{l}, ..., f^{n-1}(x)\in U_{l}$, thus
$(-\infty,0]\subseteq \bigcap\limits_{i=0}^{n-1}f^{-i}(U_{l})$.

As $f$ is a monotone increasing mapping, there exists $k > 0$ such that $f^{k}(b_{l})>a_{r}$.
We can assume $n > k > 0$.   
Consider the following two possibilities 1.1) and 1.2).

\smallskip
\smallskip

1.1) $x\in [b_{l},a_{r}]$.

It requires at most $k$ iterations so that
$f^{k}(x)\in U_{r}$. Hence,
$x\in U_{j_{0}}, f(x)\in U_{j_{1}}, ..., f^{k-1}(U_{j_{k-1}})$,
$f^{k}\in U_{r}, ..., f^{n-1}(x)\in U_{r}$ where $U_{j_{0}},U_{j_{1}}, ., U_{j_{k-1}} \in\mathcal{U}$.
Since $card(\mathcal {U})=m$, $[b_{l},a_{r}]$ can be covered by $m^{k}$ elements of
$\bigvee\limits_{i=0}^{n-1}f^{-i}(\mathcal{U})$.

\smallskip
\smallskip

1.2) $x\in (0,b_{l})$.

This is divided into  three further possibilities  as follows.
\smallskip

1.2.1) $f^{n-1}(x)>a_{r}$.

Choose $j$ with $0 < j < n$ such that $f^{j-1}(x)<b_{l}$, but $f^{j}(x) \geq b_{l}$.
Then $x\in U_{l},f(x)\in  U_{l},...,f^{j-1}(x)\in U_{l},f^{j}(x)\in U_{j_{0}},...,f^{j+k-1}(x)\in
U_{j_{k-1}},f^{j+k}(k)\in U_{r},...,f^{n-1}(x)\in
U_{r}$, where $U_{j_{0}},U_{j_{1}}, ..., U_{j_{k-1}} \in\mathcal{U}$.
Since $card(\mathcal{U})=m$,
$\bigvee\limits_{i=0}^{n-1}f^{-i}(\mathcal{U})$
has $m^{k}$ elements that cover this kind of points $x$.
\smallskip

1.2.2) $b_{l}\leq f^{n-1}(x)\leq a_{r}$.

If $f^{n-2}(x)<b_{l}$, i.e., for the last jump getting into $[b_{l},a_{r}]$, it holds
$x\in U_{l}, ..., f^{n-2}(x)\in U_{l}, f^{n-1}(x)\in U_{j_{0}}$, where
$U_{j_{0}}\in\mathcal{U}$ while $card(\mathcal{U})=m$,
there are $m$ elements of $\bigvee\limits_{i=0}^{n-1}f^{-i}(\mathcal{U})$ that
cover these kind of points $x$.

If $f^{n-3}(x)<b_{l}$ and $f^{n-2}(x)\geq b_{l}$, i.e., for the second jump from last before getting
into $[b_{l},a_{r}]$, it holds
$x\in U_{l}, ..., f^{n-3}(x)\in U_{l}, f^{n-2}(x)\in U_{j_{2}}, f^{n-1}(x)\in
U_{j_{1}}$, where $U_{j_{2}}, U_{j_{1}}\in\mathcal {U}$ while $card(\mathcal{U})=m$,
$\bigvee\limits_{i=0}^{n-1}f^{-i}(\mathcal{U})$ has $m^{2}$ elements that
cover this kind of points $x$.

Continue in this fashion ...,
if $f^{n-k}(x)<b_{l}$ and $f^{n-(k-1)}(x)\geq b_{l}$, i.e., for the $(k-1)$th jump from last
before getting into $[b_{l},a_{r}]$, it holds
$x\in U_{l}, ..., f^{n-k}(x)\in U_{l}, f^{n-(k-1)}(x)\in U_{j_{k-1}}, ., f^{n-1}(x)\in U_{j_{1}}$,
where $U_{j_{1}}, ..., U_{j_{k-1}} \in \mathcal{U}$ while $card(\mathcal{U})=m$,
$\bigvee\limits_{i=0}^{n-1}f^{-i}(\mathcal{U})$ has $m^{k-1}$ elements that
cover this kind of points $x$.

If $f^{n-(k+1)}(x)<b_{l}$ and $f^{n-k}(x)\geq b_{l}$,
i.e., jump into $[b_{l},a_{r}]$ on the $k$th,  $f^{n-1}(x)>a_{r}$ and this is Case 1.2.1).
\smallskip

1.2.3) $f^{n-1}(x)<b_{l}$.

Clearly, $x\in \bigcap\limits_{i=0}^{n-1}f^{-i}(U_{l})\in \bigvee\limits_{i=0}^{n-1}f^{-i}(\mathcal{U})$.
\smallskip
\smallskip

Hence, in Case 1 where  $0<b_{l}<a_{r}$, for any $n > k > 0$,  
it holds $N(\bigvee\limits_{i=0}^{n-1}f^{-i}(\mathcal{U}))\leq
2 + m^{k} + m^{k}+m+m^{2} + ... + m^{k-1}$, and by the definition of co-compact entropy,
$c(f,\mathcal{U})=\lim\limits_{n\rightarrow\infty}\frac
{1}{n}H(\bigvee\limits_{i=0}^{n-1}f^{-i}(\mathcal{U}))\leq
\lim\limits_{n\rightarrow\infty}\frac
{1}{n}\log(2+m^{k}+m^{k}+m+m^{2}+...+m^{k-1})=0$.

\smallskip
\smallskip

Case 2: $b_{l}<a_{r}<0$. This is similar to above Case 1.

\smallskip
\smallskip

Case 3: $b_{l}<0<a_{r}$.
For any $n > 0$  and $x\in(a_{r},+\infty)$,
we have $x\in U_{r}, f(x)\in U_{r}, ..., f^{n-1}(x) \in  U_{r}$, thus $(a_{r},+\infty)\subseteq
\bigcap\limits_{i=0}^{n-1}f^{-i}(U_{r})$.

Similarly,  for $x\in(-\infty,b_{l})$, we have
$x\in U_{l}, f(x)\in  U_{l}, ..., f^{n-1}(x)\in U_{l}$, thus  $(-\infty,b_{l})\subseteq  \bigcap\limits_{i=0}^{n-1}f^{-i}(U_{l})$.  As $\mathcal{U}$ is an open cover of $R$, there exists $U_{0}\in \mathcal{U}$ such that
$0\in U_{0},f(0)=0\in  U_{0},...,f^{n-1}(0)=0\in U_{0}$, and hence
$0\in\bigcap\limits_{i=0}^{n-1}f^{-i}(U_{0})$.

For $x\in[b_{l},a_{r}]$, $U_{0}$ as an open set of $R$ can be decomposed into a union of countably many open intervals.
Noting that $0 \in U_0$, there are two further possibilities, as given in 3.1) and 3.2) below.
\smallskip
\smallskip

3.1) The stated decomposition  of $U_0$ has an interval  $(b_{0},a_{0})$ that contains $0$, i.e., $0 \in (b_{0},a_{0})$.
Since $f$ is a monotone increasing mapping, there exists $k > 0$ such that $f^{k}(b_{0})<b_{l}$ and $f^{k}(a_{0})>a_{r}$.
Here, we can assume $n > k > 0$.   
Similar to Case 1, $[b_{l},b_{0}]$ can be covered by $m^{k}$ elements of
$\bigvee\limits_{i=0}^{n-1}f^{-i}(\mathcal{U})$,
$(b_{0},0)$ can be covered by $m^{k}+m+m^{2}+...+m^{k-1}$ elements of
$\bigvee\limits_{i=0}^{n-1}f^{-i}(\mathcal{U})$,
$(0,a_{0})$ can be covered by $m^{k}+m+m^{2}+...+m^{k-1}$ elements of
$\bigvee\limits_{i=0}^{n-1}f^{-i}(\mathcal{U})$, and
$[a_{0},a_{r}]$ can be covered by $m^{k}$ elements of
$\bigvee\limits_{i=0}^{n-1}f^{-i}(\mathcal{U})$. Hence, for any $n > k > 0$,
we have $N(\bigvee\limits_{i=0}^{n-1}f^{-i}(\mathcal{U}))\leq 3+m^{k}+m+m^{2}+...+m^{k-1}+ m^{k}+m+m^{2}+...+m^{k-1}$. Therefore, by the definition of co-compact entropy,
$c(f,\mathcal{U})=\lim\limits_{n\rightarrow\infty}\frac{1}{n}H(\bigvee\limits_{i=0}^{n-1}f^{-i}(\mathcal{U}))\leq
\lim\limits_{n\rightarrow\infty}\frac{1}{n}\log(3+m^{k}+m+m^{2}+...+m^{k-1}+ m^{k}+m+m^{2}+...+m^{k-1})=0$.
\smallskip
\smallskip

3.2) The only intervals covering $0$ are of the forms $(-\infty,a_{0})$ or $(b_{0},+\infty)$.

Consider the case $0\in(-\infty,a_{0})$. As $f$ is a monotone increasing mapping,
there exists $k>0$ such that $f^{k}(a_{0}) >a_{r}$. We can assume $n > k > 0$. 
Similar to Case 1, $(0,a_{0})$ can be covered by
$m^{k}+m+m^{2}+...+m^{k-1}$ elements of  $\bigvee\limits_{i=0}^{n-1}f^{-i}(\mathcal{U})$,
$[a_{0},a_{r}]$ can be covered by $m^{k}$ elements of
$\bigvee\limits_{i=0}^{n-1}f^{-i}(\mathcal{U})$, and it also holds
$[b_{l},0)\subseteq \bigcap\limits_{i=0}^{n-1}f^{-i}(U_{0})$. Hence, for any $n > k > 0$, we have $N(\bigvee\limits_{i=0}^{n-1}f^{-i}(\mathcal{U}))\leq 3+m^{k}+m+m^{2}+...+m^{k}$. By
the definition of co-compact entropy, we have
$c(f,\mathcal{U})=\lim\limits_{n\rightarrow\infty}\frac {1}{n}H(\bigvee\limits_{i=0}^{n-1}f^{-i}(\mathcal{U}))\leq
\lim\limits_{n\rightarrow\infty}\frac {1}{n}\log(3+m^{k}+m+m^{2}+...+m^{k} )=0$.
Therefore, when $b_{l}<0<a_{r}$,  it holds $c(f,\mathcal{U})=0$.

The case $0 \in (b_{0},+\infty)$ is similar.
\smallskip
\smallskip

Now, by Cases 1, 2 and 3, it holds $c(f,\mathcal{U})=0$.  Noting that $\mathcal{V}\prec\mathcal{U}$,
it holds $c(f,\mathcal{V})\leq c(f,\mathcal{U})=0$.
Since $\mathcal{V}$ is arbitrary, we have $c(f)=0$.

\section{\hspace*{-0.15in}. Concluding remarks}   \label{sec6}			
The investigation of dynamical systems could be tracked back to Isaac Newton's era when calculus and his
laws of motion and universal gravitation were invented or discovered, in which differential equations with time as a parameter play
a dominant role. However, it was not realized until  the end of the 19th century that the hope of solving all kinds of problems in
celestial mechanics
by following Newton's frame and methodology, e.g., the two body problem, becomes unrealistic when
Jules Henri Poincar\'{e}'s New Methods of Celestial Mechanics was publicized (shortly after this, in the early 20th century,
fundamental changes in electrodynamics occurred when Albert Einstein's historical papers  appeared:
reconciling Newtonian mechanics with 
Maxwell's electrodynamics,  separating  Newtonian  mechanics from quantum mechanics,
and extending the principle of relativity to non-uniform motion), in which the space of
all potential values of the parameters of the system is included in the analysis, and the attention to the system was
changed from individual solutions to dynamical properties of all solutions as well as the relation among all solutions. Although this approach
may not provide much information on individual solutions, it can obtain important information of most of the solutions. For example,
by taking an approach similar to that in ergodic theory, Poincar\'{e} concluded that for all Hamiltonian systems, most solutions are
stable \cite{Ye}.

The study of dynamical systems has
become a central part of mathematics and its applications since the middle of the 20th century when scientists from all
related disciplines realized the power and beauty of the geometric and qualitative techniques developed during this period
for nonlinear systems (see e.g., Robinson \cite{Robinson}).

Chaotic and random behavior of solutions of deterministic systems is now understood to be an inherent feature of
many nonlinear systems (Devaney \cite{Devaney}, 1989).
Chaos and related concepts as main concerns in mathematics and physics were investigated
through differentiable dynamical systems, differential equations, geometric structures, differential topology,
and ergodic theory etc, by e.g.,
S. Smale, J. Moser, M. Peixoto, V.I. Arnol'd, Ya. Sinai, J.E. Littlewood, M.L. Cartwright, A.N. Kolmogorov,
G.D. Birkhoff among others, and even as early
as H. Poincar\'{e} (global properties, nonperiodicity, 1900's) and J. Hadamard (stability of trajectories, 1890's).

Kolmogorov's metric entropy as an invariant of measure theoretical dynamical systems
is a powerful concept because it controls the top of the hierarchy of ergodic properties, and
plays a remarkable role in investigating the complexity and other properties of the systems. As an analogous
invariant under conjugation of topological dynamical systems,
topological entropy plays prominent  role for
the study of dynamical systems, and is often used as a measure in determining  dynamical behavior (e.g., chaos) and complexity of systems.
In particular,
topological entropy bounds measure-theoretic entropy (Goodwyn \cite{Goodwyn, Goodwyn2}). Other relations between various entropy
characterizations were extensively studied, e.g., Dinaburg \cite{Dinaburg}.
It is a common understanding that topological entropy, as a non-negative number and invariant of conjugation
in describing dynamical systems, serves a unique and unsubstitutable role in dynamics. Consequently, an appropriate
definition of topological entropy becomes important and difficult.

In theory and applications of dynamical systems, locally compact systems appear commonly, e.g., $R^n$ or other manifolds.
The introduced concept of co-compact open covers  is fundamental for describing the dynamical behaviors of systems as, for example, for locally
compact systems co-compact open sets are the neighborhoods of the infinity point in the Alexandroff  compactification and hence
admit the investigation of the dynamical properties near infinity.

The co-compact entropy introduced in this paper is defined based on the co-compact open covers.
In the special case of compact systems, this new entropy coincides with
the topological entropy introduced by Adler, Konheim and McAndrew (Sections 3 and 4).
For non-compact systems, this new entropy retains
various fundamental properties of Adler, Konheim and McAndrew's entropy (e.g., invariant under conjugation, entropy of a
subsystem does not exceed that of the whole system).

Another noticeable property of the co-compact entropy is that it is metric-independent for dynamical systems defined on metric
spaces, thus different from the entropy defined by Bowen. In particular, for the linear mapping given in Section \ref{SEC5}
(locally compact system), its co-compact entropy  is $0$, which would be at least $\log 2$ according to Bowen's definition;
as a positive entropy usually reflects certain dynamical complexity of a system, this new entropy is more appropriate.

For a dynamical system defined on a metric space, Bowen's definition may result in different entropies when
different metrics are employed. As proved in Section 4, the co-compact entropy is a lower bound for Bowen's entropies, where
the traditional Lebesgue Covering Theorem for open covers of  compact metric spaces is generalized
for co-compact open covers of non-compact metric spaces, too.

\end{document}